\documentclass[titlepage,11pt]{article}
\usepackage{latexsym}
\usepackage{amsfonts}
\usepackage{amsmath}
\usepackage{underscore}
\newcommand\blackslug{\hbox{\hskip 1pt \vrule width 4pt height 8pt depth 1.5pt
        \hskip 1pt}}
\newcommand\bbox{\hfill \quad \blackslug \bigbreak}

\def\ll{,\ldots,}

\title{Proof of the Kalai-Meshulam conjecture}

\author{Maria Chudnovsky\thanks{Supported by NSF grant DMS-1550991.
This material is based upon work supported in part by the U. S. Army
Research Laboratory and the U. S. Army Research Office under grant
number
W911NF-16-1-0404.}\\
Princeton University, Princeton, NJ 08544
\\
\\
Alex Scott\thanks{Supported by a Leverhulme Research
Fellowship.}\\
Mathematical Institute, University of Oxford, Oxford OX2 6GG, UK
\\
\\
Paul Seymour\thanks{Supported by ONR grant N00014-14-1-0084 and NSF
grant DMS-1265563.}\\
Princeton University, Princeton, NJ 08544
\\
\\
Sophie Spirkl\\
Princeton University, Princeton, NJ 08544}

\date{June 4, 2018; revised \today}

\newtheorem{thm}{}[section]

\newcommand{\Proof}{\noindent{\bf Proof.}\ \ }

\begin{document}
\maketitle
\begin{abstract} 
Let $G$ be a graph, and let $f_G$ be the sum of $(-1)^{|A|}$, over all 
stable sets $A$. If $G$ is a cycle with length divisible by three, then
$f_G= \pm 2$.
Motivated by topological considerations, G. Kalai and R. Meshulam~\cite{kalai} made the conjecture that,
if no induced cycle of a graph $G$ has length divisible by three, 
then $|f_G|\le 1$.
We prove this conjecture.
\end{abstract}

\section{Introduction}

All graphs in this paper are finite and have no loops or parallel edges.
A graph is {\em ternary} if no induced cycle has length a multiple of three;
thus, ternary graphs have no triangles. In the late 1990's 
G. Kalai and R. Meshulam~\cite{kalai}
made two conjectures about these graphs, the following (both now theorems):

\begin{thm}\label{conj0}
There exists $c$ such that every ternary graph is $c$-colourable.
\end{thm}
\begin{thm}\label{firstconj1}
For every ternary graph, the number of stable sets with even 
cardinality and the number with odd cardinality differ by at most one.
\end{thm}
The first was proved by Bonamy, Charbit and Thomass\'e~\cite{bonamy} (although 
it may be that all ternary graphs are 3-colourable, and this remains open).
A much stronger result was later proved by two of us~\cite{holeparity}:
that for all integers $p,q\ge 0$, every graph with bounded clique number and with
no induced cycle of length $p$ modulo $q$ has bounded chromatic number.

The second conjecture has remained open, and we prove it in this paper. 
We mention a few other related results (there are more in the final section):

\begin{itemize}
\item Chen and Saito~\cite{chen} proved that every non-null graph with no cycle of length divisible by three (not just induced cycles)
has a vertex of degree at most two (and so all such graphs are 3-colourable). 
\item G. Gauthier~\cite{gauthier}
found an explicit construction for all graphs with no cycle of length
divisible by three.
\item D. Kr\'al\textquoteright~asked (unpublished): is it true that in every ternary
graph with an edge, there is an edge $e$ such that the graph obtained by
deleting $e$ is also ternary? This would have implied                               
that all ternary graphs are 3-colourable, but has very recently been disproved; a counterexample was found by M. Wrochna. 
(Take the disjoint union of a 5-cycle and a 10-cycle, and join each vertex of the 5-cycle to 
two opposite vertices of the 10-cycle, in order.)
\item Kalai and Meshulam also proposed the conjecture that for all $k$ there 
exists $c$, such that, if for every induced subgraph of $G$ the number of 
even stable sets and the number of odd ones differ by at most $k$, then $G$ is $c$-colourable.
This is proved in~\cite{holeparity}.
\end{itemize}

If $G$ is a graph, and $X,Y$ are disjoint subsets of $V(G)$, let $f_G(X,Y)$
be the sum of $(-1)^{|A|}$, summed over all stable sets $A$ in $G$ that include
$X$ and are disjoint from $Y$. Our main theorem states:
\begin{thm}\label{mainthm}
If $G$ is ternary then $|f_G(\emptyset, \emptyset)|\le 1$.
\end{thm}

The proof of \ref{mainthm} is by induction on $|V(G)|$, and it follows easily that if $G$ is a minimum counterexample
then $f_G(\emptyset,\emptyset)=\pm 2$. It is very helpful to
know the value of $f_G(\emptyset,\emptyset)$, and so the proof breaks into two cases, depending whether this value 
is $2$ or $-2$. The proof for the second is obtained from the first proof by negating $f_G$ throughout, 
and we would like to say ``we may assume that $f_G(\emptyset,\emptyset)=2$ without loss of generality''; but 
this gives us a difficulty, because negating 
$f_G$ does not give a function that equals $f_H$ for some graph $H$. We overcome this as follows. 

Let $G$ be a graph, and with $f_G$ as before, let us say the functions $f_G$ and $-f_G$ are 
{\em counters} on $G$. We will prove that if $G$ is  ternary and $g$ is a counter on $G$, then 
$|g(\emptyset,\emptyset)|\le 1$. Now we are free to replace $g$ by its negative if that is convenient.

We will frequently need to talk about $g(X,Y)$ when $Y=\emptyset$; so often that it is worthwhile to make a special
convention for it. We define $g(X)=g(X,\emptyset)$ (and the same for $f_G$).

If $g$ is a counter on $G$, we say $g$ is a {\em good} counter if for 
all disjoint $X,Y\subseteq V(G)$ with $X\cup Y\ne \emptyset$:
\begin{itemize}
\item  $|g(X,Y)|\le 1$; and
\item 
$|g(X\cup \{u\},Y) -g(X\cup \{v\},Y)|\le 1$ for all $u,v\in V(G)\setminus (X\cup Y)$.
\end{itemize} 
In section \ref{sec:conj3}, we show that:
\begin{thm}\label{conj3}
If $g$ is a good counter on a graph $G$, then $|g(\{u\}) -g(\{v\})|\le 1$ for all $u,v\in V(G)$.
\end{thm}
Then in section \ref{sec:conj1}, we show that:
\begin{thm}\label{conj1}
If $g$ is a good counter on a ternary graph $G$, then $|g(\emptyset)|\le 1$.
\end{thm}
\noindent{\bf Proof of \ref{mainthm}, assuming \ref{conj3} and \ref{conj1}.\ \ }
We prove by induction on $|V(G)|$ that for every ternary graph $G$, if $g$ is a counter on $G$, then
$|g(\{u\}) -g(\{v\})|\le 1$ for all $u,v\in V(G)$, and $|g(\emptyset)|\le 1$.
Thus we may assume that these two statements hold for every proper induced subgraph of $G$. 
Now $g$ is a counter on $G$, and so $g=\pm f_G$. If the result holds for $-g$ then it holds for $g$;
so we may assume that $g=f_G$, by replacing $g$ by $-g$ if necessary.
\\
\\
(1) {\em If $X,Y\subseteq V(G)$ are disjoint, with $X\cup Y\ne \emptyset$, then
$|f_G(X,Y)|\le 1$.}
\\
\\
We may assume that $X$ is a stable set. Let 
$H$ be the graph obtained from $G$ by deleting $X\cup Y$ and deleting all vertices with a neighbour in $X$.
Thus, if $A$ is a stable set of $G$ including $X$ and disjoint from $Y$, then $A\setminus X$ is a stable set
of $H$; and conversely, if $B$ is a stable set of $H$, then $X\cup B$ is a stable set of $G$ including $X$ 
and disjoint from $Y$. 
In particular, $f_H(\emptyset) = (-1)^{|X|}f_G(X,Y)$; but from the 
inductive hypothesis, $|f_H(\emptyset)|\le 1$, and so $|f_G(X,Y)|\le 1$. This proves (1).
\\
\\
(2) {\em If $X,Y\subseteq V(G)$ are disjoint, with $X\cup Y\ne \emptyset$, and $u,v\in V(G)\setminus (X\cup Y)$, 
then 
$$|f_G(X\cup \{u\},Y) -f_G(X\cup \{v\},Y)|\le 1.$$}
We may assume that $X$ is stable.
Suppose first that $u$ has a neighbour in $X$. Then $f_G(X\cup \{u\}, Y)=0$ (because
$X\cup\{u\}$ is not a subset of any stable set). Also $|f_G(X\cup \{v\}, Y)|\le 1$, by (1), and the claim follows.
So we may assume that $u$ and similarly $v$ has no neighbour in $X$; and so $u,v\in V(H)$, if we define $H$
as before. Thus $f_G(X\cup \{u\},Y)=(-1)^{|X|}f_H(\{u\})$, and 
$f_G(X\cup \{v\},Y)=(-1)^{|X|}f_H(\{v\})$; and from the inductive hypothesis,
$|f_H(\{u\}) -f_H(\{v\})|\le 1$. It follows that
$|f_G(X\cup \{u\},Y) -f_G(X\cup \{v\},Y)|\le 1$. This proves (2).

\bigskip
From (1) and (2), $g$ is a good counter on $G$. From \ref{conj1} and \ref{conj3}, it follows that 
$|g(\{u\}) -g(\{v\})|\le 1$ for all $u,v\in V(G)$, and $|g(\emptyset)|\le 1$.
This completes the inductive proof; and \ref{mainthm} follows.~\bbox

\section{Some lemmas}
Here are a few useful lemmas. First, we observe:

\begin{thm}\label{incexc}
Let $g$ be a counter on $G$, let $X,Y\subseteq V(G)$ be disjoint, and let $Y'\subseteq Y$.
Then 
$$g(X,Y)=\sum_{Z\subseteq Y\setminus Y'}(-1)^{|Z|}g(X\cup Z, Y').$$
\end{thm}
\Proof
We may assume that $g=f_G$, by replacing $g$ by $-g$ if necessary. 
If $A$ is a stable set of $G$ including $X$ and disjoint from $Y'$, define $n_A$ to be 
$$\sum_{Z\subseteq A\cap Y} (-1)^{|A|-|Z|}.$$ 
Thus $n_A=0$ unless $A\cap Y=\emptyset$, in which case $n_A=(-1)^{|A|}$.
But 
$\sum_{Z\subseteq Y\setminus Y'}(-1)^{|Z|}f_G(X\cup Z, Y')$ is the sum of $n_A$, over all stable sets $A$ of $G$
including $X$ and disjoint from $Y'$. It follows that 
$\sum_{Z\subseteq Y\setminus Y'}(-1)^{|Z|}f_G(X\cup Z, Y')$ is the sum of $(-1)^{|A|}$ over all stable sets of $G$
that include $X$ and are disjoint from $Y$. But this sum equals $f_G(X,Y)$. This proves \ref{incexc}.~\bbox

In evaluating an expression given by \ref{incexc}, it often happens that for some number $\ell$,
$g(X\cup Z)=\ell$ for ``most'' subsets $Z\subseteq Y$, and if so the following is helpful:
\begin{thm}\label{zerosum}
Let $g$ be a counter on $G$, let $X,Y\subseteq V(G)$ be disjoint, with $Y\ne \emptyset$, and let $\ell$ be some number.
Then
$$g(X,Y)=\sum_{Z\subseteq Y}(-1)^{|Z|}(g(X\cup Z)-\ell).$$
\end{thm}
\Proof
By \ref{incexc}, 
$$g(X,Y)=\sum_{Z\subseteq Y}(-1)^{|Z|}g(X\cup Z),$$
and $\sum_{Z\subseteq Y}(-1)^{|Z|}(-\ell)=0$ since $Y\ne \emptyset$. This proves \ref{zerosum}.~\bbox

\begin{thm}\label{adjacent} 
Let $g$ be a good counter on $G$, let $X,Y\subseteq V(G)$ be disjoint, and let $v\in
V(G)\setminus (X\cup Y)$. Then $|g(X,Y)-g(X\cup \{v\},Y)|\le 1$ and $|g(X,Y)-g(X,Y\cup \{v\})|\le 1$.
\end{thm}
\Proof
We may assume that $g=f_G$. Every stable set including $X$ and disjoint from $Y$ either includes $X\cup \{v\}$
or is disjoint from $Y\cup \{v\}$, and not both. Consequently 
$$g(X,Y)=g(X\cup\{v\}, Y) + g(X,Y\cup \{v\}).$$
But $|g(Y\cup\{v\})|\le 1$ since $g$ is a good counter, 
and therefore $|g(X,Y)-g(X\cup \{v\},Y)|\le 1$; and the second claim follows similarly.~\bbox

For $X\subseteq V(G)$, let $N[X]$
denote the set of vertices in $G$ that either belong to $X$ or have a
neighbour in $X$.
We observe that
\begin{thm}\label{lonely}
Let $g$ be a counter on $G$.
If $X,Y\subseteq V(G)$ are disjoint with $g(X,Y)\ne 0$, and
$v\in V(G)\setminus (N[X]\cup Y)$, then
$v$ has a neighbour in $V(G)\setminus (N[X]\cup Y)$.
\end{thm}
\Proof We may assume that $f=f_G$, by replacing $g$ by $-g$ if necessary.
The stable sets of $G$ that include $X$ and are disjoint from $Y$ are obtained 
from the stable sets of $G\setminus (N[X]\cup Y)$ ($=H$ say) by adding the set $X$ to each such stable set; and so
$f_H(\emptyset)\ne 0$.
But $f_K(\emptyset)=0$ for every graph $K$
with a vertex of degree zero, and so $H$ has no vertex of degree zero. The result follows.~\bbox

\begin{thm}\label{Sophieconf}
Let $g$ be a good counter on $G$, let $X,Y\subseteq V(G)$ be disjoint, and let $u,v\in V(G)\setminus (X\cup Y)$.
If $g(X,Y)=g(X\cup \{u,v\}, Y)\ne 0$, then $g(X,Y)=g(X\cup \{v\}, Y)$.
\end{thm}
\Proof We proceed by induction on $|V(G)\setminus (X\cup Y)|$. By replacing $g$ by $-g$ if necessary we may
assume that $g(X,Y)>0$. 
For all disjoint $A,B\subseteq V(G)\setminus (X\cup Y)$, let $h(A,B)=g(X\cup A, Y\cup B)$ (and $h(A)$ means $h(A,\emptyset)$).
Since $g$ is a good counter it follows that $|h(\{u,v\})|\le 1$, and so 
$h(\{u,v\})=h(\emptyset)=1$. We suppose for a contradiction that $h(\{v\})\ne 1$.
Hence $u\ne v$, and $X\cup \{u,v\}$ is stable.
By \ref{adjacent}, it follows that $h(\{v\})=0$. Since $|h(\emptyset, \{u,v\})|\le 1$,
\ref{incexc} implies that 
$$h(\emptyset)-h(\{u\})-h(\{v\})+h(\{u,v\})\le 1.$$
Consequently
$h(\{u\})\ge 1$, and so $h(\{u\})= 1$. From \ref{lonely}, $v$ has a neighbour $w$.

Now $h(\emptyset,\{v\})=h(\emptyset)-h(\{v\}) = 1$, and 
$h(\{u\},\{v\})=h(\{u\})-h(\{u,v\}) = 0$, and so from the inductive hypothesis,
$h(\{u,w\},\{v\})\ne 1$. Consequently $h(\{u,w\})-h(\{u,v,w\})\ne 1$, and since
$h(\{u,v,w\})=0$, it follows that $h(\{u,w\})\ne 1$.
By \ref{adjacent}, $h(\{u,w\})=0$. 
Thus $h(\{u\},\{w\})=1$ by \ref{incexc}, since $h(\{u\})=1$. Since $h(\{v\},\{w\})=0$ and
$h(\{u,v\},\{w\})=1$ by \ref{incexc} (the first since $h(\{v,w\})=0$ and $h(\{v\})=0$,
and the second since $h(\{u,v,w\})=0$ and $h(\{u,v\})=0$), it follows from the inductive hypothesis
that $h(\emptyset,\{w\})\ne 1$, and so $h(\emptyset,\{w\})=0$ by \ref{adjacent}. Hence 
$h(\emptyset)-h(\{w\})=0$ by \ref{incexc}, and so $h(\{w\})=1$. But then
$h(\{w\},\{u\})=1$, because $h(\{u,w\})=0$; and $h(\{v\},\{u\})=-1$, since $h(\{v\})=0$
and $h(\{u,v\})=1$. This contradicts that $g$ is good, and so proves \ref{Sophieconf}.~\bbox

The next result has been independently discovered several times.
\begin{thm}\label{3col}
Let $G$ be a nonnull 
graph and let $A_1,A_2,A_3$ be the classes of a 3-colouring of $G$. Suppose that for $i = 1,2,3$, every vertex
in $A_i$ has a neighbour in $A_{i+1}$, where $A_4$ means $A_1$. Then $G$ is not ternary.
\end{thm}
\Proof
Throughout we read subscripts modulo 3.
For $i = 1,2,3$, direct each edge of $G$ between $A_i$ and $A_{i+1}$ from $A_i$ to $A_{i+1}$. Since each vertex has 
positive outdegree, the digraph we form has a directed cycle, and hence an induced directed cycle. But such a cycle
is an induced cycle of $G$, and has length a multiple of three.~\bbox

\begin{thm}\label{hypergraph}
Let $\mathcal{H}$ be a set of subsets of some set $V$, all of the same cardinality $k$; and suppose that
for every subset $X\subseteq V$ with $|X|=k+1$, if $X$ includes a member of $\mathcal{H}$ then it includes
at least two such members. Then there is a partition $P_1\ll P_n$ of $V$ with $P_1\ll P_n$ all nonempty,
such that
for all distinct $u,v\in V$, either there exists $i\in \{1\ll n\}$ with $u,v\in P_i$, or there exists
$B\in \mathcal{H}$ with $u,v\in B$, and not both.
\end{thm}
\Proof 
Say two vertices $u,v\in V$ are {\em equivalent} if either $u=v$, or:
\begin{itemize}
\item there is no member of $\mathcal{H}$ containing both $u,v$; and
\item for each $C\subseteq V\setminus \{u,v\}$, $C\cup \{u\}\in \mathcal{H}$ if and only if $C\cup \{v\}\in \mathcal{H}$.
\end{itemize}
We claim that this is an equivalence relation. To see this, we may assume that $u,v,w\in V(G)$ are distinct, and $v$ is equivalent
to both $u$ and $w$; and we must show that $u,w$ are equivalent. If there exists $B\in \mathcal{H}$ containing $u,w$, then $v\notin B$
(since $u,v$ are equivalent) and so $(B\setminus \{u\}\cup \{v\}\in \mathcal{H}$ (since $(B\setminus \{u\})\cup \{u\}\in \mathcal{H}$
and $u,v$ are equivalent), and so this is a member of $\mathcal{H}$ containing $v,w$, a contradiction. Thus there is no such $B$.
Let $C\subseteq V\setminus \{u,w\}$, with $C\cup \{u\}\in \mathcal{H}$. Consequently $v\notin C$, and $C\cup \{v\}\in \mathcal{H}$
(because $u,v$ are equivalent), and consequently $C\cup \{w\}\in \mathcal{H}$ (since $v,w$ are equivalent). Similarly
$C\cup \{u\}\in \mathcal{H}$ if and only if $C\cup \{v\}\in \mathcal{H}$. This proves that equivalence is indeed an equivalence
relation.

We claim that for
all distinct $u,v\in V$, if they do not belong to the same equivalence class then some member of $\mathcal{H}$
contains both $u,v$. To see this, since $u,v$ are not equivalent, if no member of $\mathcal{H}$
contains both $u$ and $v$, then we may assume (exchanging $u,v$ if necessary)
that there exists $C\subseteq V\setminus \{u,v\}$
such that $C\cup \{u\}\in \mathcal{H}$ and $C\cup \{v\}\notin \mathcal{H}$.
Thus $|C|=k-1$, and since $C\cup \{u,v\}$ includes a member of $\mathcal{H}$, by hypothesis it includes at least two members.
But since no member of $\mathcal{H}$ contains both $u,v$, and $C\cup \{v\}\notin \mathcal{H}$, this is impossible.
This proves \ref{hypergraph}.~\bbox

\section{The value on distinct vertices}\label{sec:conj3}

In this section we prove \ref{conj3}. Thus, throughout this section, let $g$ be a good counter on a graph $G$. 
For $i = -1,0,1$ let $A_i$
be the set of vertices $v$ of $G$ such that $g(\{v\})=i$. Thus $A_{-1}, A_0, A_1$ are disjoint and have
union $V(G)$. We need to show that one of $A_{-1}$, $A_1$ is empty, and so we assume for a contradiction
that they are both nonempty. We will prove a series of statements about $G,g$.
We begin with:

\begin{thm}\label{stable}
The following hold:
\begin{itemize}
\item $g(\emptyset)=0$;
\item $G$ is connected;
\item $A_1, A_{-1}$ are both stable sets;
\item there is not both an edge between $A_1,A_0$ and an edge between $A_{-1}, A_0$.
\end{itemize}
\end{thm}
\Proof
Since there exists $v\in A_1$, and hence with $g(\{v\})=1$, we deduce from 
\ref{adjacent} 
that $g(\emptyset)\ge 0$, and similarly $g(\emptyset)\le 0$. This proves the first statement.

For the second statement,
we may assume (replacing $g$ by $-g$ if necessary) that $g=f_G$. 
By assumption, there exist $u_i\in V(G)$ with $g(\{u_i\})=i$, for $i = 1,-1$.
Suppose that $G$ is not connected, and 
let $G_1$ be a component of $G$ containing $u_1$, and let $G_2$ be obtained from $G$ by deleting $G_1$.
Write $g_i$ for $f_{G_i} (i=1,2)$. Thus for disjoint $X,Y\subseteq V(G)$, 
$$g(X,Y)=g_1(X\cap V(G_1),Y\cap V(G_1))g_2(X\cap V(G_2),Y\cap V(G_2)).$$
In particular, since $0=g(\emptyset)=g_1(\emptyset)g_2(\emptyset)$,
one of $g_1(\emptyset),g_2(\emptyset)$ is zero.

Since 
$g(\{u_1\})=g_1(\{u_1\})g_2(\emptyset)$, it follows that $g_2(\emptyset)\ne 0$,
and so $g_1(\emptyset)=0$. In particular, $G_1$ is the unique component $C$ of $G$ such that
$f_C(\emptyset)=0$, and so $u_{-1}\in V(G_1)$. Thus 
$g(\{u_{-1}\})=g_1(\{u_{-1}\})g_2(\emptyset)$,
and so one of $g_1(\{u_{1}\}), g_1(\{u_{-1}\})$ equals $1$ and the other equals $-1$,
contradicting that $g$ is good. This proves the second statement.

For the third, suppose that $u,v\in A_1$ are adjacent. By \ref{incexc},
$$g(\emptyset, \{u,v\})=g(\emptyset)-g(\{u\})-g(\{v\})+g(\{u,v\});$$
but the last term is zero since $u,v$ are adjacent, and since $u,v\in A_1$ and $g(\emptyset)=0$,
we deduce that 
$g(\emptyset, \{u,v\})=2$, contradicting that $g$ is good.

For the fourth statement, suppose that $u_1\in A_1$ is adjacent to $v_1\in A_0$, and $u_{-1}\in A_{-1}$
is adjacent to $v_{-1}\in A_0$. Suppose first that $g(\{v_1,u_{-1}\})=0$. Then by two applications of \ref{incexc},
$g(\{u_{-1}\},\{v_1\})= g(\{u_{-1}\})-g(\{u_{-1},v_1\})=-1$, and
$g(\{u_1\},\{v_1\})= g(\{u_1\})-g(\{u_1,v_1\})=1$ (since $u_1, v_1$ are adjacent),
contradicting that $g$ is good. This proves that $g(\{v_1,u_{-1}\})\ne 0$, and so 
$g(\{v_1,u_{-1}\})=-1$ by \ref{adjacent}. Similarly $g(\{v_{-1},u_1\})=1$ (and in particular,
$v_1\ne v_{-1}$). But by \ref{incexc}, 
$$g(\{v_1\},\{u_1,u_{-1}\}) = g(\{v_1\})-g(\{v_1,u_1\})
-g(\{v_1,u_{-1}\})+g(\{v_1,u_1,u_{-1}\});$$
and since $g(\{v_1\})=0$ and $g(\{v_1,u_1\})=g(\{v_1,u_1,u_{-1}\})=0$ (since $u_1,v_1$
are adjacent) it follows that $g(\{v_1\},\{u_1,u_{-1}\})=1$. Similarly $g(\{v_{-1}\},\{u_1,u_{-1}\})=-1$,
contradicting that $g$ is good.
This proves \ref{stable}.~\bbox

In the same notation, because of the third statement of \ref{stable}, we may assume (replacing $g$ by $-g$ if necessary)
that there are no edges between $A_{-1}$ and $A_0$. Let $B_1$ be the set of vertices $v\in A_0$ such that 
$g(\{u,v\})=1$ for each $u\in A_1$ and $g(\{u,v\})=0$ for each $u\in A_{-1}$; and let $B_2$
be the set of vertices $v\in A_0$ such that  
$g(\{u,v\})=0$ for each $u\in A_1$ and $g(\{u,v\})=-1$ for each $u\in A_{-1}$.

\begin{thm}\label{partition}
Every vertex in $A_0$ belongs to one of $B_1, B_2$.
\end{thm}
\Proof Let $v\in A_0$, and for $i=1,-1$ let $u_i\in A_i$. Not both $g(\{v,u_1\})=1$ and
$g(\{v,u_{-1}))=-1$, since $g$ is good. Suppose that neither of these holds. Then $g(\{v,u_1\})=0$
and $g(\{v,u_{-1}))=0$, by \ref{adjacent}. Then by two applications of \ref{incexc},
$g(\{u_1\},\{v\})=g(\{u_1\})-g(\{u_1,v\}) = 1$, and 
$g(\{u_{-1}\},\{v\})=g(\{u_{-1}\})-g(\{u_{-1},v\}) = -1$,
contradicting that $g$ is good. It follows that either $g(\{v,u_1\})=1$ and $g(\{v,u_{-1}\})=0$,
or $g(\{v,u_1\})=0$ and $g(\{v,u_{-1}\})=-1$. Since this holds for all $u_1, u_{-1}$,
it follows that $v\in B_1\cup B_2$. This proves \ref{partition}.~\bbox

\begin{thm}\label{partition2}
$A_0$ is empty.
\end{thm}
\Proof Suppose that $B_1$ is nonempty. Since $G$ is connected by \ref{stable}, and there are no edges between $B_1$
and $A_1$ (since $g(\{u,v\})=1$ for each $u\in A_1$ and $v\in B_1$, from the definition of $B_1$), it follows
that there is an edge between $B_1, B_2$; say between $b_1\in B_1$ and $b_2\in B_2$. Also, since there are no edges
between $A_{-1}$ and $A_0$, and $A_{-1}\ne \emptyset$, there is an edge between $A_{-1}$ and $A_1$; say between
$a_{-1}\in A_{-1}$ and $a_1\in A_1$. Then $g(\{a_1,b_1\})=1$, $g(\{a_1,b_2\})=0$,
$g(\{a_{-1},b_1\})=0$, and $g(\{a_{-1},b_2\})=-1$. Since $b_1$ is adjacent to $b_2$,
\ref{incexc} implies that 
$$g(\{b_1\}, \{a_1,a_{-1}\})=g(\{b_1\})-g(\{b_1,a_1\})-g(\{a_1,a_{-1}\})+
g(\{b_1,a_1,a_{-1}\})=-1,$$ and similarly $g(\{b_2\}, \{a_1,a_{-1}\})=1$, contradicting that $g$ is good. 
This proves that $B_1=\emptyset$.

Now suppose that $B_2\ne \emptyset$. Since $G$ is connected, there is an edge between $B_2,A_1$, say between $v\in B_2$
and $u_1\in A_1$. Choose $u_{-1}\in A_{-1}$. By three applications of \ref{incexc},
$g(\emptyset,\{u_1\})=g(\emptyset)-g(\{u_1\}) = -1$,
$g(\{v\},\{u_1\})=g(\{v\})-g(\{v,u_1\}) =0$, and
$g(\{v,u_{-1}\},\{u_1\})=g(\{v,u_{-1}\})-g(\{v,u_1,u_{-1}\}) =-1$,
contrary to \ref{Sophieconf}. Thus $B_2=\emptyset$ and so $A_0=\emptyset$. This proves \ref{partition2}.~\bbox

Now we prove \ref{conj3}, which we restate:
\begin{thm}\label{conj3again}
If $g$ is a good counter on a graph $G$, then $|g(\{u\}) -g(\{v\})|\le 1$ for all $u,v\in V(G)$.
\end{thm}
\Proof 
In the previous notation,  \ref{partition2} and \ref{stable} imply that $G$ is bipartite, and 
$(A_1, A_{-1})$ is a bipartition. We recall that $g(\emptyset) = 0$.
\\
\\
(1) {\em Every vertex of $G$ has degree at least two.}
\\
\\
Since $G$ is connected by \ref{stable}, all vertices have degree at least one; suppose that $v\in A_1$ has only one neighbour
$u\in A_{-1}$ say. Since $G$ is connected and $|V(G)|\ge 3$, $u$ has another neighbour $v'\in A_1$. Now 
$g(\{v'\})=1$, and since $v\in V(G)\setminus N[\{v'\}]$, \ref{lonely} implies that 
$v$ has a neighbour in $V(G)\setminus N[\{v'\}]$, a contradiction. This proves (1).
\\
\\
(2) {\em There is a subset $X\subseteq A_1$ with $g(X)=0$.}
\\
\\
Choose $v\in A_1$, and let $X=A_1\setminus \{v\}$. Since $v\in V(G)\setminus N[X]$, and $v$ has no neighbour in 
$V(G)\setminus N[X]$ (by (1)), \ref{lonely} implies that $g(X)=0$. This proves (2).

\bigskip

For $i = 1,-1$ let $k_i>0$ be minimum such that some subset $B$ of $A_i$ with cardinality $k_i$ satisfies $g(B)\ne i$.
Thus $k_i\ge 2$; and by \ref{adjacent}, $g(B)=0$ or $i$ for each subset $B\subseteq A_i$ with $|B|=k_i$.
\\
\\
(3) {\em For $i = 1,-1$, $k_i$ is odd.}
\\
\\
Since $g$ is good, $|g(\emptyset, B_i)|\le 1$; and so by \ref{zerosum},
$$\left|\sum_{Z\subseteq B_i}(-1)^{|Z|}(g(Z)-i)\right|\le 1.$$
But $g(Z)=i$ for all $Z\subseteq B_i$ with $Z\ne B_i,\emptyset$, and zero if $Z= B_i,\emptyset$; and 
consequently $|-i-i(-1)^{k_i}|\le 1$, and so
$k_i$ is odd. This proves (3).

\bigskip

Let $\mathcal{H}_i$ be the set of all subsets $B$ of $A_i$ such that $|B|=k_i$ and $g(B)=0$.
\\
\\
(4) {\em For every subset $X$ of $A_i$ with cardinality $k_i+1$, if $X$ includes a member of $\mathcal{H}_i$
then it includes at least two such members.}
\\
\\
Let $X=\{v_0\ll v_{k_i}\}$, and suppose that $\{v_1\ll v_{k_i}\}$ is the only member of $\mathcal{H}_i$ included in $X$.
Then $g(X)\ne i$, by \ref{Sophieconf}, and $g(X)\ne -i$ by \ref{adjacent}; so $g(X)=0$.
Let $Y=\{v_2\ll v_{k_i}\}$.
By \ref{zerosum} and (3):
\begin{eqnarray*}
g(\emptyset, Y)&=&\sum_{Z\subseteq Y}(-1)^{|Z|}(g(Z)-i)=-i,\\
g(\{v_0\},Y)&=&\sum_{Z\subseteq Y}(-1)^{|Z|}(g(Z\cup \{v_0\})-i)=0,\\
g(\{v_0,v_1\},Y)&=&\sum_{Z\subseteq Y}(-1)^{|Z|}(g(Z\cup \{v_0,v_1\})-i)= -(-1)^{|Y|}i=-i,
\end{eqnarray*}
contrary to \ref{Sophieconf}. This proves (4).
\\
\\
(5) {\em There exist $B_i\in \mathcal{H}_i$ for $i=1,-1$, such that 
there are two edges of $G$ between $B_1$ and $B_{-1}$ with no end in common.}
\\
\\
By (4) and \ref{hypergraph}, there is a partition $P_1\ll P_m$ of $A_1$ such that every two vertices in $A_1$ either belong to
the same $P_i$ or to some member of $\mathcal{H}_1$, and not both; and let $Q_1\ll Q_n\subseteq A_{-1}$ be 
defined analogously.
Say $P_i, Q_j$ are adjacent if there is an edge in $G$ between a vertex in $P_i$
and a vertex in $Q_j$. Since $m,n\ge 2$ and each $P_i$ is adjacent to some $Q_j$ and vice versa, there are
distinct $P_1,P_2$ (say) and distinct $Q_1,Q_2$ such that $P_1$ is adjacent to $Q_1$ and $P_2$ to $Q_2$.
Choose $p_i\in P_i$ and $q_i\in Q_i (i=1,2)$ such that $p_1q_1$ and $p_2q_2$ are edges of $G$. Since $p_1,p_2$
do not belong to the same one of $P_1\ll P_m$, there exists $B_1\in \mathcal{H}_1$ containing $p_1,p_2$; and similarly
there exists $B_{-1}\in \mathcal{H}_{-1}$
containing $q_1,q_2$. This proves (5).

\bigskip
For $i = 1,-1$ choose $B_i$ as in (5).
\\
\\
(6) {\em For $i=1,-1$, let $X_i\subseteq B_i$ with $\emptyset\ne X_i\ne B_i$. Then $g(X_1\cup X_{-1})=0$.}
\\
\\
Suppose not, and for $i = 1,-1$ choose $X_i\subseteq B_i$ with $\emptyset\ne X_i\ne B_i$,  
with $X_1\cup X_{-1}$ minimal such that $g(X_1\cup X_2)\ne 0$.
We may assume that $g(X_1\cup X_{-1})=1$, by replacing $g$ by $-g$
if necessary. By \ref{incexc} and the minimality of $X_1\cup X_{-1}$,
$$g(X_1,X_{-1}) = g(X_1)+(-1)^{|X_{-1}|}g(X_1,X_{-1})=1+(-1)^{|X_{-1}|},$$
and so $|X_{-1}|$ is odd; and similarly $|X_1|$ is even. Choose $u\in X_1$ and $v\in X_{-1}$.
Then by three applications of \ref{incexc},
\begin{eqnarray*}
g(X_1\setminus\{u\}, X_{-1}\setminus \{v\})&=&g(X_1\setminus\{u\})=1,\\
g((X_1\cup \{v\})\setminus\{u\}, X_{-1}\setminus \{v\})&=&0,\\
g(X_1\cup \{v\}, X_{-1}\setminus \{v\})&=&(-1)^{|X_{-1}\setminus \{v\}|}g(X_1\cup X_2)=1,
\end{eqnarray*}
contrary to \ref{Sophieconf}. This proves (6).

\bigskip

Choose $C_1\subseteq B_1$ maximal such that either $C_1=\emptyset$ or $g(C_1\cup B_{-1})\ne 0$,
and choose $C_{-1}\subseteq B_{-1}$ maximal such that 
either $C_{-1}=\emptyset$ or $g(C_{-1}\cup B_{1})\ne 0$. It follows that $|C_i|\le k_i-2$ for $i = 1,-1$,
since there is a 2-edge matching between $B_1, B_{-1}$.
For $i=1,-1$ let $D_i=B_i\setminus C_i$, and let $C=C_1\cup C_{-1}$ and $D=D_1\cup D_{-1}$.
\\
\\
(7) {\em If $C_1\ne \emptyset$ then $g(C_1\cup B_{-1})=1$; and if $C_{-1}\ne \emptyset$ then
$g(C_{-1}\cup B_{1})=-1$.}
\\
\\
Since $g(C_1, B_{-1})\ne 2$ (because $g$ is good), and $g(C_1\cup Z)=0$ for all $Z\subseteq B_{-1}$
with $Z\ne \emptyset, B_{-1}$ by (6), \ref{incexc} implies that 
$g(C_1)+(-1)^{|k_{-1}|}g(C_1\cup B_{-1})\le 1$. But $g(C_1)=1$ (since $C_1\ne \emptyset$),
and $k_1$ is odd, and so $g(C_1\cup B_{-1})=1$. Similarly if $C_{-1}\ne \emptyset$
then $g(C_{-1}\cup B_{1})=-1$. This proves (7).
\\
\\
(8) {\em One of $C_1,C_{-1}$ is empty.}
\\
\\
Suppose they are both nonempty. 
By \ref{incexc},
$$g(C, D)= \sum_{Z\subseteq D}(-1)^{|Z|}g(C\cup Z).$$
But for $Z\subseteq D$, $g(C\cup Z)\ne 0$ only if $Z$ includes one of $D_1,D_{-1}$ by (6), and only if 
one of $Z\cap B_1, Z\cap B_{-1}$ is empty (from the definition of $C_1, C_{-1}$); that is, only if $Z$ is one
of $D_1,D_{-1}$. These two sets are distinct, since they are nonempty. Consequently
$$g(C,D)= (-1)^{|D_1|}g(B_{1}\cup C_{-1}) + (-1)^{|D_{-1}|}g(B_{-1}\cup C_{1})$$
and so by (7), $g(C,D)= (-1)^{|D_1|+1} + (-1)^{|D_{-1}|}$. Since $|g(C,D)|\le 1$ (because $g$ is good) it follows that
$|D_1|,|D_{-1}|$ have the same parity.

Choose $u\in D_1$ and $v\in D_{-1}$. Then by \ref{incexc},
$$g(C\cup \{u\}, D\setminus \{u,v\})= \sum_{Z\subseteq D\setminus \{u,v\}}(-1)^{|Z|}g(C\cup \{u\}\cup Z).$$
But for $Z\subseteq D\setminus \{u,v\}$, $g(C\cup \{u\}\cup Z)\ne 0$ only if $Z=D_1\setminus \{u\}$ (by (6)
and the definition of $C_1, C_{-1}$) and so 
$$g(C\cup \{u\}, D\setminus \{u,v\})=(-1)^{|D_1\setminus \{u\}|}g(B_1\cup C_{-1})= (-1)^{|D_1|}.$$
Similarly 
$$g(C\cup \{v\}, D\setminus \{u,v\})=(-1)^{|D_2\setminus \{v\}|}g(B_{-1}\cup C_{1})= (-1)^{|D_{-1}|+1}.$$
Since $|D_1|,|D_{-1}|$ have the same parity, 
one of $g(C\cup \{u\}, D\setminus \{u,v\}), g(C\cup \{v\}, D\setminus \{u,v\})$ equals $1$ and the other equals $-1$,
contradicting that $g$ is good. This proves (8).

\bigskip

From (8) we may assume that $C_{-1}=\emptyset$ (replacing $g$ by $-g$ if necessary).
\\
\\
(9) {\em $|D_1|$ is odd.}
\\
\\
To prove this, we may assume that $C_1\ne \emptyset$.
By \ref{incexc},
$$g(C_1, B_{-1}\cup D_1) = \sum_{Z\subseteq B_{-1}\cup D_1}(-1)^{|Z|}g(C_1\cup Z).$$
But, by (6), for $Z\subseteq B_{-1}\cup D_1$, $g(C_1\cup Z)$ is nonzero only if $Z\subseteq D_1$ or $Z=B_{-1}$; and then it has value
1 if $Z\subseteq D_1$ and $Z\ne D_1$;
0 if  $Z=D_1$; and 
1 if $Z=B_{-1}$.
Thus $g(C_1, B_{-1}\cup D_1)=(-1)^{|D_1|+1}+(-1)^{|B_{-1}|}$ and since $|B_{-1}|$ is odd by (5), and 
$|g(C_1, B_{-1}\cup D_1)|\le 1$ since $g$ is good, it follows that $|D_1|$ is odd. This proves (9).

\bigskip

Now $|C_1|\le |B_1|-2$ as we saw. Choose $u\in D_1$ and $v\in B_{-1}$, and let $W=(D_1 \cup B_{-1})\setminus \{u,v\}$.
By \ref{incexc}, 
$$g(C_1\cup \{u\}, W)= \sum_{Z\subseteq W}(-1)^{|Z|}g(C_1\cup \{u\}\cup Z).$$
But for $Z\subseteq W$, $g(C_1\cup \{u\}\cup Z)$ is nonzero only if
$Z\subseteq D_1$, and in that case it has value $1$ if $Z\ne D_1\setminus \{u\}$, and $0$ if $Z=D_1\setminus \{u\}$.
Since $|D_1|\ge 2$, it follows that
$$g(C_1\cup \{u\}, W)=(-1)^{|D_1|} =-1$$ 
since $|D_1|$ is odd by (9).
On the other hand, by \ref{incexc},
$$g(C_1\cup \{v\}, W)= \sum_{Z\subseteq W}(-1)^{|Z|}g(C_1\cup \{v\}\cup Z).$$
We claim that $g(C_1\cup \{v\}, W)=1$. To see this there are two cases, depending whether $C_1\ne \emptyset$ or not.
First, suppose that $C_1\ne \emptyset$. 
Then for $Z\subseteq W$, $g(C_1\cup \{v\}\cup Z)$ is nonzero only if
$Z=B_{-1}\setminus \{v\}$, by (6) and the maximality of $C_1$; so 
$$g(C_1\cup \{v\}, W)= (-1)^{|B_1|-1}g(C_1\cup B_{-1})= 1,$$
by (7) and (3),
contradicting that $g$ is good. 
Now suppose that $C_1=\emptyset$. Then, again by (6), for $Z\subseteq W$, $g(C_1\cup \{v\}\cup Z)$ is nonzero only if
$Z$ is a proper subset of $B_{-1}\setminus \{v\}$, and in that case it has value $-1$. Consequently
$$g(C_1\cup \{v\}, W)= (-1)^{|B_{-1}\setminus \{v\}|}g(C_1\cup B_{-1})=1,$$
again contradicting that $g$ is good.
This proves \ref{conj3again}.~\bbox

\section{The value on the null set}\label{sec:conj1}

In this section we prove \ref{conj1}, thereby completing the inductive proof of \ref{mainthm}. 
We need to show that if $g$ is a good counter on a ternary graph $G$, then $|g(\emptyset)|\le 1$.
The proof is divided into several steps. 
We may assume the statement is false, for a contradiction; and by replacing $g$ by $-g$ if necessary,
we may assume that $g(\emptyset)\ge 2$. Throughout this section, $G$ is a counterexample to \ref{conj1},
and $g$ is a good counter on $G$, with $g(\emptyset)\ge 2$.
\begin{thm}\label{step1}
The following hold:
\begin{itemize}
\item $g(\emptyset)=2$;
\item $g(\{v\})=1$ for every vertex $v\in V(G)$; and
\item $G$ is connected.
\end{itemize}
\end{thm}
\Proof
Let $v\in V(G)$; since $g$ is good, it follows that $|g(\{v\})|\le 1$, and so \ref{adjacent} implies that
$g(\{v\})=1$ and $g(\emptyset)=2$. This proves the first two statements.

Suppose that $G$ is not connected, let $G_1$ be a component of $G$ and let $G_2$ be obtained from $G$ by deleting $V(G_1)$.
Since $f_{G_1}(\emptyset)= \pm g(\emptyset, V(G_2))$, and $g$ is good, it follows that
$|f_{G_1}(\emptyset)|\le 1$, and similarly $|f_{G_2}(\emptyset)|\le 1$.
But
$$g(\emptyset) = \pm f_G(\emptyset)=\pm f_{G_1}(\emptyset)f_{G_2}(\emptyset),$$
a contradiction. This proves the third statement, and so proves \ref{step1}.~\bbox

In particular, 
if $u,v\in V(G)$ are distinct, then since $g(\{u\})=1$ by the second statement of \ref{step1}, 
it follows that $g(\{u,v\})\in \{0,1\}$ by \ref{adjacent}.
Let $H$ be the graph with vertex set $V(G)$ in which distinct $u,v$ are adjacent if $g(\{u,v\})=1$.
\begin{thm}\label{step2}
Every component of $H$ is a complete graph, and $H$ has at least two and at most four components.
\end{thm}
\Proof
Suppose the first statement is false.
Then there are three distinct vertices $u,v,w\in  V(H)$ such that $uv,vw\in E(H)$ and $uw\notin E(H)$.
From \ref{adjacent}, $g(\{u,w\})=0$. Now 
\begin{eqnarray*}
g(\emptyset,\{w\})&=&g(\emptyset)-g(\{w\}) = 1,\\
g(\{v\},\{w\})&=&g(\{v\})-g(\{v,w\}) = 0\\
g(\{u,v\},\{w\}) &=& g(\{u, v\})-g(\{u,v,w\});
\end{eqnarray*}
and by \ref{Sophieconf}, $g(\{u,v\},\{w\})\ne 1$. Consequently $g(\{u,v,w\})=1$.
But then $g(\{w\}) = 1$, $g(\{u,w\})=0$ and $g(\{u,v,w\})=1$, contrary to \ref{Sophieconf}.
This proves that every component of $H$ is a complete graph.

Since each edge of $H$ joins two vertices that are nonadjacent in $G$, it follows that $H$ has at least two components.
Suppose it has at least five. Since $G$ is connected, there is a vertex of $H$ that has neighbours (in $G$) in at least 
two components of $H$. Thus we can choose $v_1\ll v_5\in V(G)$, all in different components of $H$, where 
$v_1$ is adjacent (in $G$) to $v_2, v_3$.  Let $a,b,c\in \{v_1\ll v_5\}$ be distinct. Since $|g(\emptyset,\{a,b,c\})|\le 1$,
and $g(\{a,b\})=0$ (because $a,b$ belong to different components of $H$, and by \ref{adjacent}),
and the same for $\{a,c\}$ and $\{b,c\}$, it follows from \ref{incexc} that
$|2-3+0-g(\{a,b,c\})|\le 1$, and so $g(\{a,b,c\})\ne 1$. Hence $g(\{a,b,c\})\in \{0,-1\}$ for every triple $a,b,c$ of distinct members of
$\{v_1\ll v_5\}$. 

Note that since $v_1v_2, v_1v_3\in E(G)$, it follows that $g(\{v_1,v_2,v_i\})=0$ for
every $i \{3,4,5\}$ and $g(\{v_1,v_3,v_j\})=0$ for every $j \in \{2,4,5\}$.
Let $\mathcal{T}$ be the set of all subsets $T\subseteq \{v_1\ll v_5\}$ with $|T|=3$ and $g(T)=-1$.
Thus  $g(T)=0$ for all triples $T\notin \mathcal{T}$.
Since 
$|g(\emptyset,\{v_1,v_2,v_3,v_4\})|\le 1$, it follows from \ref{incexc} that 
$\{v_2,v_3,v_4\}\in \mathcal{T}$, and similarly $\{v_2,v_3,v_5\}\in \mathcal{T}$.

Suppose that $\{v_1,v_4,v_5\}\notin \mathcal{T}$.  
Now \ref{incexc} implies that
$$g(\emptyset, \{v_1,v_2,v_4,v_5\})= 2-4+0-g(\{v_2,v_4,v_5\}),$$ 
and so
$\{v_2,v_4,v_5\}\in \mathcal{T}$, and similarly $\{v_3,v_4,v_5\}\in \mathcal{T}$. But then
$$g(\{v_5\},\{v_2,v_3,v_4\})=-2-g(v_2,v_3,v_4,v_5)\le -1$$ and 
$g(\{v_1\},\{v_2,v_3,v_4\})=1$, contradicting that $g$ is good.
Thus $\{v_1,v_4,v_5\}\in \mathcal{T}$. 

If also $\{v_2,v_4,v_5\}\in \mathcal{T}$ then $g(\{v_4,v_5\},\{v_1,v_2\})=2$, contradicting that 
$g$ is good; so $\{v_2,v_4,v_5\}\notin \mathcal{T}$, and similarly $\{v_3,v_4,v_5\}\notin \mathcal{T}$.
Since $g(\{v_2,v_3\},\{v_4,v_5\})\le 1$, it follows that $g(\{v_2,v_3,v_4,v_5\})=-1$. But then
$g(\{v_4\},\{v_2\}) = 1$, $g(\{v_4,v_5\},\{v_2\})=0$ and $g(\{v_3,v_4,v_5\},\{v_2\})=1$, contrary to \ref{Sophieconf}.
This proves \ref{step2}.~\bbox

\begin{thm}\label{crossset}
Let $C_1,C_2$ be distinct components of $H$, and let $X\subseteq C_1\cup C_2$.
Suppose that 
\begin{itemize}
\item $X\cap C_1, X\cap C_2\ne \emptyset$;
\item $g(X)\ne 0$; and
\item for all $X'\subseteq X$, if $g(X')\ne 0$ then either $X'=X$ or $X'\subseteq C_1$ or $X'\subseteq C_2$.
\end{itemize}
If $|X\cap C_1|>1$ then there is a subset $B\subseteq X\cap C_1$ with $g(B)=0$.
\end{thm}
\Proof Let $X_i=X\cap C_i$ for $i = 1,2$; and suppose there is no $B\subseteq X_1$ with $g(B)=0$. 
From \ref{adjacent} it follows that $g(B)=1$ for all nonempty subsets
$B$ of $X_1$, and in particular, $g(X_1)=1$. Let $g(X)=i=\pm 1$. 
Because of the third bullet of the hypothesis, \ref{incexc} implies that
$$g(X_1,X_2) = \sum_{Z\subseteq X_2}(-1)^{|Z|}g(X_1\cup Z)=g(X_1)+(-1)^{|X_2|}i;$$
and since $g(X_1,X_2)\le 1$, it follows that $(-1)^{|X_2|}i=-1$, that is, $|X_2|$ is odd if $i=1$, and even if $i=-1$.
Choose $u\in X_1$ and $v\in X_2$; then by \ref{incexc}, $g(X_1\setminus\{u\}, X_2\setminus \{v\}) = 1$ (since $|X_1|>1$),
$g(X_1\cup \{v\}\setminus\{u\}, X_2\setminus \{v\}) = 0$, and by \ref{incexc},
$$g(X_1\cup \{v\}, X_2\setminus \{v\}) =\sum_{Z\subseteq X_2\setminus \{v\}}(-1)^{|Z|}g(X_1\cup Z\setminus \{v\})
=(-1)^{|X_2|-1}g(X)=1,$$
contrary to \ref{Sophieconf}. This proves \ref{crossset}.~\bbox

Let $C$ be a component of $H$, and let $D\subseteq C$. We say that $B\subseteq D$ is a {\em base}
of $D$ if $g(B)\ne 1$ and there is no $B'\subseteq D$ with $|B'|<|B|$ and with $g(B')\ne 1$.
\begin{thm}\label{bases}
Let $C$ be a component of $H$, and let $D\subseteq C$. 
\begin{itemize}
\item If there is a vertex $v$ of $G$ such that all its neighbours belong to $D$, then $D$ has a base.
\item If $B$ is a base of $D$ then $g(B)=0$, and $|B|$ is even and at least four.
\item If $D$ has a base, of cardinality $k$ say, then every subset of $D$ of cardinality $k+1$
includes two bases of $D$, and so every vertex of $D$ belongs to a  base of $D$.
\item If $D$ has a base, of cardinality $k$, then there is a partition of $D$ into at least $k$ nonempty sets
$D_1\ll D_n$, such that for all distinct $u,v\in D$, there is a base of $D$ containing both $u,v$ if and only if
$u,v$ do not belong to the same set $D_i$.
\end{itemize}
\end{thm}
\Proof  For the first statement, suppose that all neighbours of $v$ belong to $D$. If $V(G)=C\cup \{v\}$,
then $v$ is adjacent to all other vertices (since no vertex has degree zero, by \ref{lonely}), contradicting that
$g(\emptyset)=2$. Thus we may choose $u\notin C\cup \{v\}$. By \ref{lonely}, $g(\{u\}, D)=0$,
but $g(\{u\})=1$, and so by \ref{incexc}, there exists a nonempty subset $Z\subseteq D$ such that $g(Z\cup \{u\})\ne 0$.
Since $C$ is the vertex set of a component of $H$, it follows that $|Z|\ge 2$.
From \ref{crossset}, there exists $B\subseteq Z$ with $g(B)=0$. This proves the first statement.

For the second, let $B$ be a base of $D$. Then $g(B)\ne 1$ by hypothesis, and in particular $|B|\ge 3$,
since $B\subseteq C$. For every $B'\subseteq B$ with $B'\ne \emptyset, B$, we have $g(B')=1$, and since
there is such a choice of $B'$ with $|B'|=|B|-1$, \ref{adjacent} implies that $g(B)\ne -1$; and hence
$g(B)=0$ since $g$ is good. But by \ref{zerosum},
$$g(\emptyset,B)=\sum_{Z\subseteq B}(-1)^{|Z|}(g(Z)-1)=(g(\emptyset)-1)+(-1)^{|B|}(g(B)-1)
=1-(-1)^{|B|},$$
and so $|B|$ is even. This proves  the second statement.

For the third, let $B$ be a base of $D$, with $|B|=k$ say; it suffices to prove that for all $v\in D\setminus B$,
$B\cup \{v\}$ includes at least two bases of $D$. Let $X=B\cup \{v\}$, and choose $u\in B$. 
Thus $g(X\setminus \{u,v\})=1$
and $g(X\setminus \{v\})=0$, so by \ref{Sophieconf}, $g(X)\ne 1$. We may assume that
$g(X\setminus \{u\})=1$, and so by \ref{adjacent}, $g(X)=0$.
By \ref{incexc}, $g(\emptyset, X\setminus \{u,v\})=1$ and $g(\{u,v\}, X\setminus \{u,v\})=1$, so by
\ref{Sophieconf}, $g(\{v\}, X\setminus \{u,v\})=1$. Hence by \ref{incexc}, since $|X|\ge 3$, there exists 
$Z\subseteq X\setminus \{u,v\}$ with $g(Z\cup \{v\})\ne 1$. Then $|Z|\le |B|$, and since
$B$ is a base for $D$, it follows that $Z$ is minimal with $g(Z)\ne 1$, and hence $Z$ is another
base for $D$. This proves the third statement.

The fourth statement follows from \ref{hypergraph}. This proves \ref{bases}.~\bbox

We call a partition $D_1\ll D_n$ as in the fourth statement of \ref{bases} the {\em induced partition} of $D$,
and the sets $D_1\ll D_n$ are called its {\em classes}.
(If the partition exists then it is unique, as is easily seen.)

\begin{thm}\label{matching}
Let $C_1,C_2$ be distinct components of $H$, and for $i = 1,2$,  let $D_i\subseteq C_i$, including a base for $D_i$.
Then for one of $i = 1,2$, there is a class of the induced partition of $D_i$ that meets all edges between $D_1$ and $D_2$.
\end{thm}
\Proof
Let the induced partition of $D_1$ have classes $P_1\ll P_m$, and let the induced partition of $D_2$
have classes $Q_1\ll Q_n$. We may assume that there is no $i\in \{1\ll m\}$ such that all edges between
$D_1,D_2$ have an end in $P_i$, and there is no $j\in \{1\ll n\}$ similarly. By K\"onig's theorem, there exist
distinct $i_1,i_2\in \{1\ll m\}$ and distinct $j_1,j_2\in \{1\ll n\}$ such that there is an edge between $P_{i_1}$
and $Q_{j_1}$, and an edge between $P_{i_2}$ and $Q_{j_2}$. Hence there is a base $B_1$ for $D_1$ and a base $B_2$
for $D_2$, such that there are two edges of $G$ between $B_1,B_2$ with no end in common.
\\
\\
(1) {\em Suppose that there exists $M_1\subseteq B_1$ with $g(B_2\cup M_1)\ne 0$, and choose $M_1$
maximal with this property. Then $|M_1|\le |B_1|-2$, and $g(B_2\cup M_1)=-1$, and $|M_1|$ is odd.}
\\
\\
Since there are two edges of $G$ between $B_1,B_2$ with no end in common, and both have an end in $B_2$,
it follows that neither has an end in $M_1$, and so $|M_1|\le |B_1|-2$. Let $A_1=B_1\setminus M_1$.
By \ref{incexc},
$$g(M_1,B_2)=\sum_{Z\subseteq B_2}(-1)^{|Z|}g(M_1\cup Z).$$
But for $Z\subseteq B_2$, $g(M_1\cup Z)\ne 0$ only if $Z=\emptyset$ or $Z=B_2$, by \ref{crossset}.
Consequently $g(M_1,B_2)=g(M_1) + (-1)^{|B_2|}g(M_1\cup B_2)$.
But $g(M_1)=1$ and $|B_2|$ is even, so $g(M_1\cup B_2)= -1$ since $g$ is good.
Now by \ref{incexc},
$$g(M_1,A_1\cup B_2)=\sum_{Z\subseteq A_1\cup B_2}(-1)^{|Z|}g(M_1\cup Z).$$
But for $Z\subseteq A_1\cup B_2$, $g(M_1\cup Z)\ne 0$ only if $Z\subseteq A_1$
or $Z=B_2$; and so 
$$g(M_1,A_1\cup B_2)=(-1)^{|A_1|}(g(B_1)-1)+ (-1)^{|B_2|}g(M_1\cup B_2).$$
Since $|B_2|$ is even, $g(B_1)=0$ and $g(M_1\cup B_2)=-1$, it follows that
$g(M_1,A_1\cup B_2)=(-1)^{|A_1|+1}-1$, and so $|A_1|$ is odd, and therefore so is $|M_1|$.
This proves (1).
\\
\\
(2) {\em There do not exist $M_1\subseteq B_1$ and $M_2\subseteq B_2$ with 
$g(B_2\cup M_1),  g(B_1\cup M_2)\ne 0$ and with $M_1,M_2$ both nonempty.}
\\
\\
Suppose such sets $M_1,M_2$ exist and choose them maximal. Let $A_i=B_i\setminus M_i$ for $i = 1,2$.
By (1), $g(B_2\cup M_1),  g(B_1\cup M_2)=-1$, and $|M_1|, |M_2|$ are odd.
Thus $|A_1|$ and $|A_2|$ are odd, and so $g(M_1\cup M_2, A_1\cup A_2)= 2$ by \ref{incexc}, a contradiction, 
This proves (2).
\\
\\
(3) {\em $g(X)=0$ for all $X\subseteq B_1\cup B_2$ with $X\cap B_1, X\cap B_2$ both nonempty.}
\\
\\
Suppose not; then from \ref{crossset}, and by exchanging $C_1,C_2$ if necessary, we may assume that 
there exists $M_1\subseteq B_1$, nonempty, with $g(B_2\cup M_1)\ne 0$. Choose $M_1$ maximal.
By (1), $g(B_2\cup M_1)=-1$ and $|M_1|$ is odd. Let $A_1=B_1\setminus M_1$, and choose $u\in A_1$.
Choose $v\in B_2$. Then by \ref{incexc}, since $A_1\setminus \{u\}\ne \emptyset$, it follows that
$g(M_1\cup \{u\}, (A_1\cup B_2)\setminus \{u,v\})=-1$ and $g(M_1\cup \{v\}, (A_1\cup B_2)\setminus \{u,v\})=1$,
contradicting that $g$ is good. This proves (3).

\bigskip

From (3), \ref{incexc} implies that $g(\emptyset, B_1\cup B_2)=-2$, a contradiction. This proves \ref{matching}.~\bbox

\begin{thm}\label{cloneclass}
Let $C_1,C_2$ be distinct cliques of $H$, and suppose there is a base for $C_2$. Let $D_1\ll D_n$
be the induced partition of $C_2$. Then there is no $i\in \{D_1\ll D_n\}$ such that every edge of $G$ between
$C_2$ and $V(G)\setminus (C_1\cup C_2)$ has an end in $D_i$.
\end{thm}
\Proof
Suppose there is such a value of $i$, say $i=1$. Let $A_1$ be the set of vertices in $C_1$ with neighbours in $C_2$.
Now $n\ge 4$ (by the last statement of \ref{bases}); choose $v\in D_2$. Thus all neighbours of $v$ belong to $C_1$, and hence to $A_1$. By the first statement 
of \ref{bases}, there is a base for $A_1$. By \ref{matching}, there is a set $X$ that meets all edges between 
$A_1$ and $C_2$, and $X$ is either a class of the induced partition of $C_2$ or a class of the induced partition of $A_1$.
The first is impossible since there are at least four classes of the induced partition of $C_2$, and each such class 
different from $D_1$ meets an edge
between $C_2$ and $A_1$ (because it meets some edge, and it has no edge to $V(G)\setminus (C_1\cup C_2)$ 
from the choice of $D_1$). Also the second is impossible, since each class of the induced partition of $A_1$ has an
edge to $C_2$, from the definition of $A_1$. This proves \ref{cloneclass}.~\bbox

Now we complete the proof of \ref{conj1}, which we restate:
\begin{thm}\label{conj1again}
If $g$ is a good counter on a ternary graph $G$, then $|g(\emptyset)|\le 1$.
\end{thm}
\Proof
In the same notation as before, we know that $H$ has two, three or four components. Suppose it has only two, say $C_1,C_2$.
By the first statement of \ref{bases}, there are bases for $C_1$ and for $C_2$, contrary to \ref{cloneclass}.

Now suppose that 
$H$ has exactly three components $C_1,C_2,C_3$. By \ref{3col} we may assume that some vertex $v\in C_2$
has no neighbour in $C_1$, and so by \ref{bases}, there is a base for $C_3$. Suppose that there is also a base
for $C_2$. By \ref{matching}, by exchanging $C_2,C_3$ if necessary, we may assume that there is a class of the
induced partition of $C_2$ that meets all edges between $C_2,C_3$, contrary to \ref{cloneclass}. Thus, neither
of $C_1,C_2$ have bases. By \ref{bases}, every vertex in $C_1\cup C_2$ has a neighbour in $C_3$. We recall that
$v\in C_2$
has no neighbour in $C_1$. Since $C_1$ has no base, it follows that $g(C_1)=1$, and so by \ref{lonely},
$v$ has a neighbour, $u$ say, with no neighbour in $C_1$. But then all neighbours of $u$ are in $C_2$, and so by \ref{bases},
there is a base for $C_2$, a contradiction.

Finally, suppose that $H$ has four components $C_1\ll C_4$. Let $K$ be the graph with vertex set $\{1\ll 4\}$
in which distinct $i,j$ are adjacent if there is an edge of $G$ between $C_i,C_j$. Since $G$ is connected, it follows that
every vertex of $K$ has nonzero degree. 
Suppose that $K$ has a 2-edge matching; then by renumbering $C_1\ll C_4$ we may assume that there 
exist $u_i\in C_i$ for $1\le i\le 4$
such that $u_1u_2, u_3u_4\in E(G)$.
But then $g(\emptyset,\{u_1,u_2,u_3,u_4\})=-2$ by \ref{incexc}, a contradiction.
Thus $K$ has no 2-edge matching, and since every vertex of $K$ has nonzero degree, we may assume that every edge of $K$
is incident with $1$, and so all edges of $G$ have an end in $C_1$.

For $i = 2,3,4$, let $X_i$ be the set of vertices in $C_1$ with no neighbour
in $C_i$.
By the first 
statement of \ref{bases}, there is a base for $C_1$. By \ref{cloneclass}, there is no base for $C_2$, 
and similarly none for $C_3, C_4$; and so by the first statement of \ref{bases}, every vertex of $C_1$ has neighbours in 
at least two of $C_2, C_3, C_4$. In particular, for all distinct $i,j\in \{2,3,4\}$ every vertex in $X_i$ has a neighbour in $C_j$.

Since $g(C_2)\ne 0$, \ref{lonely} implies that 
for all distinct 
$i,j\in \{2,3,4\}$, every vertex in $C_i$ has a neighbour in $X_j$. Make a digraph $J$ 
with vertex set $C_2\cup C_3\cup C_4$
in which for $i=2,3,4$ and $u\in C_i$ and $v\in C_{i+1}$ (where $C_5$ means $C_2$), there is an edge of $J$ from $u$ to $v$
if $u,v$ has a common neighbour in $X_{i-1}$ (where $X_1$ means $X_4$). Every vertex has positive outdegree in $J$,
and so $J$ has an induced directed cycle. Let $K$ be such a cycle, with vertices (in order):
$$a_1,b_1,c_1,a_2,b_2,c_2\ll a_k, b_k, c_k, a_1$$
where $a_1\ll a_k\in C_2$, $b_2\ll b_k\in C_3$ and $c_1\ll c_k\in C_4$.
For each $i$ with $1\le i\le k$, there exists $x_i\in X_4$ adjacent in $G$ to $a_i, b_i$, and $y_i\in X_2$ adjacent to 
$b_i,c_i$, and $z_i\in X_3$ adjacent to $c_i, a_{i+1}$ (where $a_{k+1}$ means $a_1$). Also, for each such $i$, $x_i$
has no other neighbours in $V(K)$; it is nonadjacent to each $a_j$ because $x_i\in X_4$,
and nonadjacent to the remaining vertices of $V(K)$ since $K$ is induced. A similar statement holds for the $y_i$'s 
and $z_i$'s. Consequently the subgraph of $G$ induced on
$$\{a_i,b_i,c_i,x_i,y_i,z_i:1\le i\le k\}$$
is an induced cycle of length $6k$, contradicting that $G$ is ternary. This proves that $H$ does not have four components,
and so proves \ref{conj1again} and hence \ref{mainthm}.~\bbox

\section{Connections}

In the 1990s, Kalai and Meshulam made a number of interesting conjectures concerning Betti numbers of graphs 
(see \cite{kalai}), and two of us proved some of
them in an earlier paper~\cite{holeparity}. But another one connects with the work of this paper.
The {\em independence complex $I(G)$} of a graph $G$ is the simplicial complex whose faces are the stable sets of vertices of $G$. 
Let $b_i$ denote the $i$th Betti number of $I(G)$ and let $b(G)$ denote the sum of the Betti numbers.   
The Kalai-Meshulam conjecture that concerns us here is:

\begin{thm}  
{\bf Conjecture: }For a graph $G$,  we have that $|b(H)| \le 1$ for every induced subgraph $H$, 
if and only if $G$ has no induced cycle of 
length divisible by 3.
\end{thm}

If $|b(H)| \le 1$ for every induced subgraph $H$, then $G$ has no induced cycle of
length divisible by 3, since $b(H)=2$ for every cycle $H$ of length divisible by three. For the converse, suppose
$G$ has no such induced cycle. Then by \ref{mainthm}, $|f_G(\emptyset)|\le 1$, but we need to prove that $b(G)\le 1$.
Now $f_G(\emptyset)$ is the Euler characteristic
of $I(G)$, and in particular there is a connection between $f_G(\emptyset)$ and $b(G)$. It is a basic theorem from homology theory
that the Euler characteristic of $I(G)$ is the alternating sum of the Betti numbers of $I(G)$ 
(see~\cite{hatcher}).
It follows that $|f_G(\emptyset)|\le b(G)$; but this inequality is in the wrong direction for us, and the conjecture remains open.

The difference between the numbers of odd and even stable sets has also appeared in statistical physics.
Let us define the polynomial
$$I_G(z)=\sum_{I}z^{|I|},$$
where the sum is over stable sets $I$ in $G$.  This polynomial is known in combinatorics as  
the {\em independent set polynomial} and statistical physics as the {\em partition function of the hard-core lattice gas} (see, for instance, \cite{AS}).  

We see that $I_G(-1)$ is the number of even stable sets minus the number of odd stable sets.    
The question of when $|I_G(-1)|\le 1$  has been the focus of considerable study, particularly on the square lattice 
(see \cite{BLN, FSvE, J}).

\end{document}